\documentclass[twoside,12pt]{article}
\usepackage{amsmath,amstext,amsgen,amsbsy,amsopn}
\usepackage{amssymb,amsfonts}
\usepackage{amsthm}
\usepackage{latexsym,amsxtra,euscript,amscd}

\usepackage{float}
\usepackage{xcolor}

\usepackage{titlesec}
\titleformat*{\section}{\normalfont\Large\bfseries\color{red}}
\titleformat*{\subsection}{\normalfont\large\bfseries\color{red}}
\titleformat*{\subsubsection}{\normalfont\normalsize\bfseries\color{red}}

\newtheoremstyle{mystyle}{3pt}{3pt}
{\itshape\color{black}}{}
{\bfseries\color{blue}}{.}{.5em}{}
\theoremstyle{mystyle}

\usepackage[font={color=brown,bf}]{caption}

\usepackage[colorlinks=true,urlcolor=blue,linkcolor=black]{hyperref}

\usepackage{graphics,graphicx}
\graphicspath{{../Fig/},{Fig/}}

\usepackage{enumitem}

\newtheorem{theorem}{Theorem}
\newtheorem{lemma}{Lemma}[section]

\newtheorem{remark}{Remark}

\newcommand{\Fl}[1]{\left\lfloor{\frac{{#1}}{2}}\right\rfloor+1}
\newcommand{\FL}[2]{\left\lfloor{\frac{{#1}}{2}}\right\rfloor#2}
\newcommand{\Cl}[1]{\left\lceil{\frac{{#1}}{2}}\right\rceil}

\usepackage{booktabs}
\newcommand{\ra}[1]{\renewcommand{\arraystretch}{#1}}

\newcommand*{\ADRnl}{University of Sopron,  Institute of Mathematics, Hungary. \texttt{nemeth.laszlo@uni-sopron.hu}}
\newcommand*{\ADRszl}{University J. Selye, Department of Mathematics, Kom\'arno, Slovakia  and University of Sopron,  Institute of Mathematics, Hungary. \texttt{szalay.laszlo@uni-sopron.hu}}

\newcommand{\seqnum}[1]{\href{https://oeis.org/#1}{\rm \underline{#1}}}

\newcommand*{\TIT}{\color{blue}Sequences involving square zig-zag shapes}
\title{\bf \TIT}

\author{L\'aszl\'o N\'emeth\footnote{\ADRnl}, L\'aszl\'o Szalay\footnote{\ADRszl} \thanks{Supported by Hungarian National Foundation for Scientific Research Grant No.~128088.}}
\date{\today}

\begin{document}
	 
\maketitle \thispagestyle{empty}

\begin{abstract}
	We define a so-called square $k$-zig-zag shape as a part of the regular square grid. Considering the shape as a $k$-zig-zag digraph, we give values of its vertices according to the number of the shortest paths from a base vertex. It provides several  integer sequences, whose higher-order homogeneous recurrences are determined by the help of a
special matrix recurrence. \\[1mm]
 	The  final  publication  is  available  at \href{https://cs.uwaterloo.ca/journals/JIS/vol24.html}{Journal of Integer Sequences}.
\end{abstract}


\section{Introduction}

The present paper studies diagonal and zig-zag paths on a particular $k+1$ wide, infinite part of the usual square lattice, and along these paths we determine linear recurrence sequences that are mostly defined in the 
{\it On-Line Encyclopedia of Integer Sequences} (OEIS, \cite{OEIS}) without combinatorial interpretations. In this manner our investigation, among others, gives them geometrical and combinatorial background. The consideration of zig-zag shapes is not an isolated challenge. For example, Baryshnikov and Romik \cite{Barys} examined 
the so-called Young diagrams, which are similar to our construction, and defined a kind of `zig-zag' numbers by the help of the alternating permutations. 
Stanley \cite{Stanley} published a survey in which he dealt with the `zig-zag' shapes and the alternating permutations. 
Recently, Ahmad et al.~\cite{Ahmad} studied  some graph-theoretic properties of special zig-zag polyomino chains.

This article proceeds our work on the variants of Pascal's arithmetic triangle for which there are several approaches to generalize or extend (see, for instance, \cite{BSz}). One recent variation is called hyperbolic Pascal triangles, see \cite{BNSz}, where we described how to construct them as we follow and generalize the connection between the classical Pascal triangle and the Euclidean regular square mosaic. The principal idea behind is to leave the Euclidean plane because in the hyperbolic plane there exist infinitely many regular mosaics. These are characterized by Schl\"afli's symbol $\{p,q\}$, where $p$ and $q$ denote two positive integers ($p,q\ge3$). The parameters indicate that for regular $p$-gons of cardinality exactly $q$ meet at each vertex. A well-defined part of each regular mosaic leads to an infinite graph such that each vertex possesses a value, giving the number of distinct shortest paths from the fixed base vertex. The structured version of the graph with labeled vertices is called hyperbolic Pascal triangle. 

More generally, if we consider an infinite connected graph $\cal G$,
the number of the shortest paths from a fixed base vertex to the other
vertices analogously generates its `${\cal G}$-Pascal triangle'. In this
paper, we present a natural family of graphs called zig-zag shape, and
investigate its properties, especially the `${\cal G}$-Pascal triangle'
of the graph. In this sense, we give a non-trivial example of `${\cal
	G}$-Pascal triangles', which may be the object of further studies.

Returning to the hyperbolic Pascal triangle linked to mosaic $\{4,5\}$ one finds the Fibonacci sequence  (\seqnum{A000045} in the OEIS \cite{OEIS}) by going along a specific zig-zag path. 
The authors  proved in \cite{BNSz,NSz_recu} that all the integer linear homogeneous recurrence sequences $\{f_i\}_{i\geq0}$ defined by
\begin{equation*}
	f_{i}=\alpha f_{i-1}\pm f_{i-2},\qquad (i\geq2),
\end{equation*} 
\noindent where $\alpha\in \mathbb{N}$, $\alpha \geq2$,  and  $f_0<f_1$ are positive integers with $\gcd(f_0,f_1)=1$,  appear along corresponding zig-zag paths in this hyperbolic Pascal triangle. 
This interesting result also inspired us to examine zig-zag paths on certain parts of the Euclidean square mosaic.	

\section{Square zig-zag shapes}\label{sec:2}

\begin{figure}[H]
	\centering
	\includegraphics[scale=1]{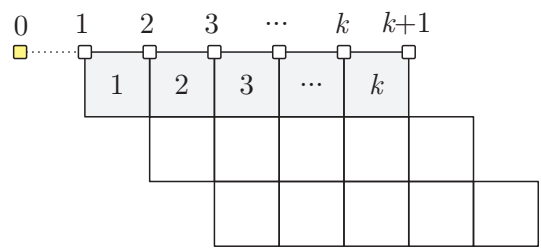} 
	\caption{Zig-zag shape}
	\label{fig:zig-zag_shape}
\end{figure}

\begin{figure}[H]
	\centering
	\includegraphics[scale=0.50]{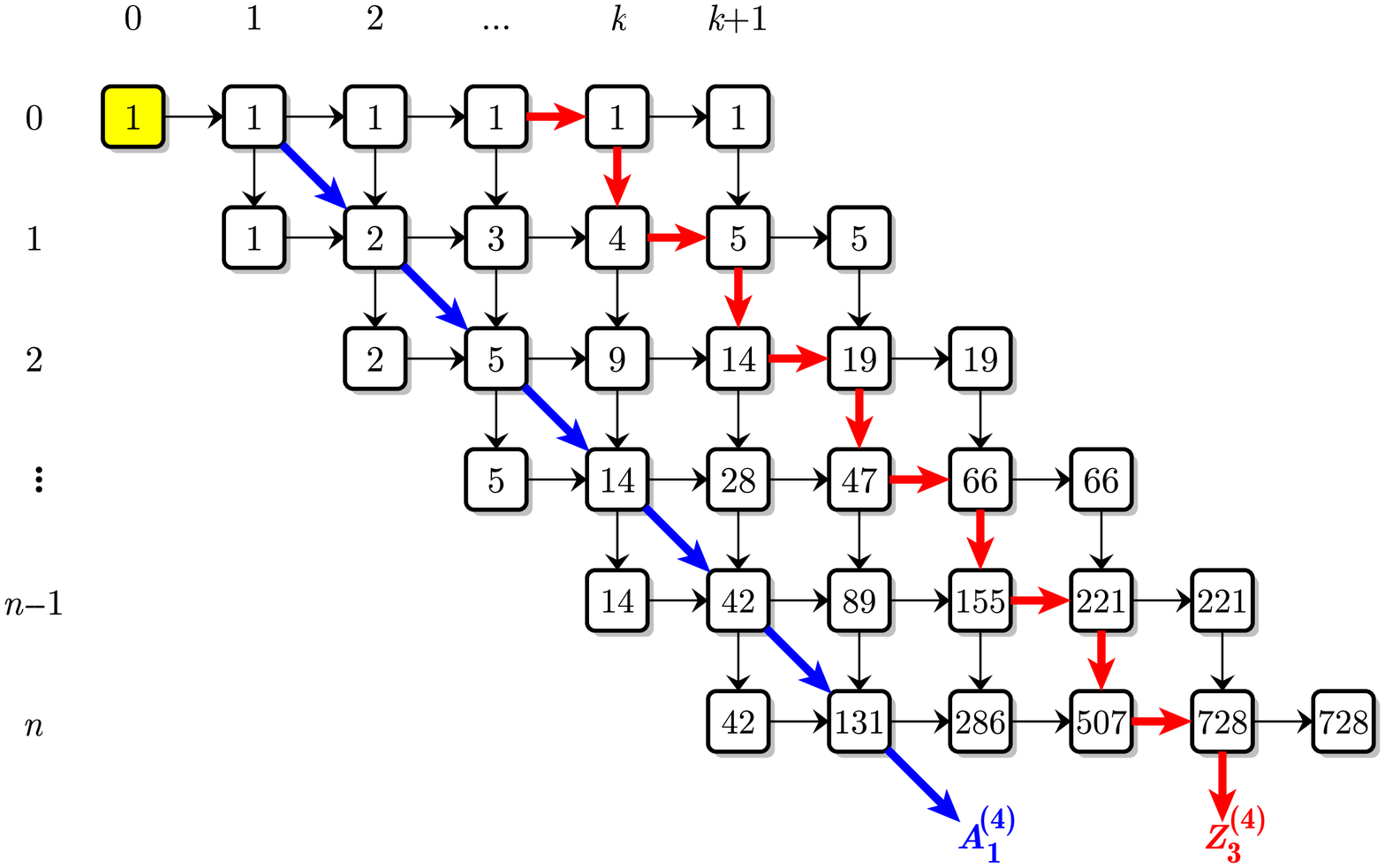} 
	\caption{Square 4--zig-zag digraph ($k=4$)}
	\label{fig:zigzag_k4}
\end{figure}

Consider the Euclidean square lattice and take $k$ consecutive pieces of squares. This is the $0$th layer of the $k$--zig-zag shape. The upper corners are the $1$st, $2$nd, $\ldots$, $k$th  and $(k+1)$st vertices according to Figure~\ref{fig:zig-zag_shape}. Extend this by an extra $0$th vertex, which is the base vertex. We color it by yellow in the figures, and we join it to the $1$st vertex by an extra edge. We denote the vertices of the $0$th line by small boxes in Figure~\ref{fig:zig-zag_shape}. Now move the $0$the layer to reach the right-down position in the square lattice to obtain the $1$st layer, and repeat this procedure with the latest layer infinitely many times. Thus, we define the square $k$--zig-zag shape or graph, where $k\geq1$ is the size of the array. Finally, we label the vertices such that a label gives the number of different shortest paths from the base vertex.  Figure~\ref{fig:zigzag_k4} illustrates the first few layers of the square 4--zig-zag digraph, the vertices are denoted by shaded boxes with their label values and the directed edges are the black arrows. (Certain black arrows are re-colored by red for some reason. There are also particular blue arrows in the Figure; their role will be discussed later.)      
Let $a_{i,j}$ denote the label of the vertice located in $i$th row and $j$th position  ($0\le j\le k+1$, $0\le i$). Clearly, the fundamental rule of the construction is given by
\begin{equation}\label{eq:base}
	a_{i,j}=\begin{cases}
		1, & \text{ if } i=0; \\ 
		a_{i-1,1}, & \text{ if } j=0, 1\leq i;\\ 
		a_{i,j-1}+a_{i-1,j+1}, & \text{ if } 1\leq j\leq k, 1\leq i; \\ 
		a_{i,k}, & \text{ if } j=k+1, 1\leq i. 
	\end{cases}
\end{equation}

For fixed $k\geq1$ and given $0\leq j\leq k+1$, let $A_j^{(k)}$ be the sequence defined by $A_j^{(k)}=(a_{i,j})_{i=0}^{\infty}$. The sequence $A_j^{(k)}$ is the $j$th right-down diagonal sequence of the square $k$--zig-zag shape. 
In Figure~\ref{fig:zigzag_k4}, the blue arrows represent the sequence $A_1^{(4)}$. 
We found  $A_0^{(k)}= (1,A_{1}^{(k)})$ and $A_k^{(k)}= A_{k+1}^{(k)}$.

Let $Z_j^{(k)}$, $j\in \{0,1,\ldots,k\}$ be the $j$th
zig-zag sequence of the  square $k$--zig-zag shape,
where $Z_j^{(k)}$ is the merged sequence of $A_j^{(k)}$ and
$A_{j+1}^{(k)}$. (In Figure~\ref{fig:zigzag_k4}, the red arrows
represent the zig-zag sequence $Z_3^{(4)}$.) More precisely,
$Z_j^{(k)}=(z_{i,j})_{i=0}^{\infty}$, where
\begin{equation}\label{eq:defbij}
	z_{i,j}=\begin{cases}
		a_{\ell,j}, & \text{ if } i=2\ell; \\ 
		a_{\ell,j+1}, & \text{ if } i=2\ell+1. 
	\end{cases}
\end{equation}
Since  $Z_0^{(k)}$ and $Z_k^{(k)}$ are  the `double' of  $A_0^{(k)}$ and $A_k^{(k)}$, respectively, usually we examine sequences for $j\in \{1,2,\ldots,k-1\}$. 

Now we record the two main theorems of this paper. The second one is a simple corollary of the first one. 
\begin{theorem}[Main 1]\label{th:main1}
	Given $k\geq1$. Then all the right-down diagonal sequences $A_j^{(k)}$ for  $j\in \{0,1,\ldots,k,k+1\}$ have the same $(\Fl{k})$-th order homogeneous linear recurrence relation 
	\begin{equation*}
		a_{n,j}= \sum_{i=0}^{\FL{k}{}}  (-1)^{i} \binom{k+1-i}{i+1} a_{n-1-i,j}, \qquad n\geq \Fl{k}.
	\end{equation*}
\end{theorem}

\begin{theorem}[Main 2]\label{th:main2}
	Fixing $k\geq1$, the zig-zag sequences $Z_j^{(k)}$ for  $j\in \{0,1,\ldots,k\}$ satisfy a $(2\FL{k}{+2})$-th order homogeneous linear recurrence relation given by
	\begin{equation*}
		z_{n,j}= \sum_{i=0}^{\FL{k}{}}  (-1)^{i} \binom{k+1-i}{i+1} z_{n-1-2i,j}, \qquad n\geq 2\FL{k}{+2}.
	\end{equation*}
\end{theorem}

For example, in case $k=4$ (see Figure~\ref{fig:zigzag_k4}) the sequences
\begin{equation*}
	\begin{aligned}
		A^{(4)}_0 &= (1,1, 2, 5, 14, 42, 131, 417, 1341, 4334, 14041, \ldots) &=&\  \text{\seqnum{A080937}},\\ 
		A^{(4)}_1 &= (1, 2, 5, 14, 42, 131, 417, 1341, 4334, 14041, \ldots) &=&\   \text{\seqnum{A080937}}\ (i\geq 1), \\ 
		A^{(4)}_2 &= (1, 3, 9, 28, 89, 286, 924, 2993, 9707, 31501, \ldots) &=&\   \text{\seqnum{A094790}},\\ 
		A^{(4)}_3 &= (1, 4, 14, 47, 155, 507, 1652, 5373, 17460, 56714,  \ldots) &=&\  \text{\seqnum{A094789}}, \\ 
		A^{(4)}_4 &= (1, 5, 19, 66, 221, 728, 2380, 7753, 25213, 81927, \ldots) &=&\   \text{\seqnum{A005021}},
	\end{aligned}
\end{equation*}
possess the common recurrence relation
\begin{equation*}
	a_{n,j}= 5 a_{n-1,j}-6 a_{n-2,j}+a_{n-3,j}, \qquad n\geq 3,
\end{equation*}
and 
\begin{equation*}
	\begin{aligned}
		Z^{(4)}_1 &= (1, 1, 2, 3, 5, 9, 14, 28, 42, 89, 131, 286, 417, \ldots) &=& \text{ not in the OEIS}, \\ 
		Z^{(4)}_2 &= (1, 1, 3, 4, 9, 14, 28, 47, 89, 155, 286, 507, 924, \ldots) &=&\ \text{\seqnum{A006053}}\ (i\geq 3), \\ 
		Z^{(4)}_3 &= (1, 1, 4, 5, 14, 19, 47, 66, 155, 221, 507, 728, 1652, \ldots)&=& \text{ not in the OEIS},
	\end{aligned}
\end{equation*}
all satisfy the recursive rule 
\begin{equation}
	z_{n,j}= 5 z_{n-2,j}-6 z_{n-4,j}+z_{n-6,j}, \qquad n\geq 6. \label{eq:znj}
\end{equation}

For more examples see Section~\ref{sec:example}. 

\begin{remark}
	The orders of recurrence relations in Theorem~\ref{th:main1} and \ref{th:main2}  are not always minimal.  Some sequences could have lower order, for example, the minimal recurrence relation of the sequence $Z^{(4)}_2$ is $z_{n,j}=  z_{n-1,j}+ 2 z_{n-2,j}-z_{n-3,j}$. Of course, its 6-order extension is \eqref{eq:znj}.
\end{remark}

\section{Recurrence relations of the square zig-zag shapes}

\begin{figure}[H]
	\centering
	\includegraphics[scale=0.6]{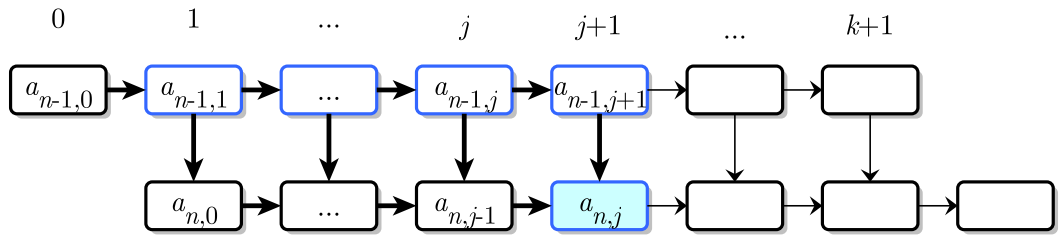} 
	\caption{Two consecutive rows}
	\label{fig:zigzag_a_ij}
\end{figure}

First look at Figure~\ref{fig:zigzag_a_ij} and the fundamental rule \eqref{eq:base}. We find that any item $a_{n,j}$, $(n\geq1)$ is the sum of the certain items  of  $(n-1)$st row.
More precisely, if $0<j<k+1$, then 
\begin{equation}\label{eq:partsum}
	a_{n,j}=a_{n-1,j+1}+a_{n,j-1}=a_{n-1,j+1}+a_{n-1,j}+a_{n,j-2}= \cdots= \sum_{\ell=1}^{j+1}a_{n-1,\ell}.\end{equation}

Consider (\ref{eq:partsum}) for all $j\in\{1,2,\dots, k+1\}$. We obtain the system  
\begin{equation}\label{eqsystem}
	\begin{aligned}
		a_{n,1} &= a_{n-1,1}+a_{n-1,2}\\ 
		a_{n,2}  &= a_{n-1,1}+a_{n-1,2}+a_{n-1,3}\\ 
		a_{n,3}  &= a_{n-1,1}+a_{n-1,2}+a_{n-1,3}+a_{n-1,4}\\ 
		&\phantom{=}\quad   \vdots \\
		a_{n,k-1}  &= a_{n-1,1}+a_{n-1,2}+a_{n-1,3}+a_{n-1,4}+\cdots+a_{n-1,k-1}+a_{n-1,k}\\ 
		a_{n,k}  &= a_{n-1,1}+a_{n-1,2}+a_{n-1,3}+a_{n-1,4}+\cdots+a_{n-1,k-1}+a_{n-1,k}+a_{n-1,k+1}\\
		a_{n,k+1}  &= a_{n-1,1}+a_{n-1,2}+a_{n-1,3}+a_{n-1,4}+\cdots+a_{n-1,k-1}+a_{n-1,k}+a_{n-1,k+1},
	\end{aligned}
\end{equation}
which includes homogenous linear recurrence sequences.

In matrix form 
\begin{equation*}
	\mathbf{v}_{n}=\mathbf{M}\cdot \mathbf{v}_{n-1}, \quad n\geq 1,
\end{equation*}
where 
$$
\mathbf{v}_{n}=
\begin{pmatrix}
	a_{n,1} \\ 
	a_{n,2} \\ 
	a_{n,3}\\ 
	\vdots\\ 
	a_{n,k} \\ 
	a_{n,k+1}
\end{pmatrix}
, 
\;\;
\mathbf{v}_{0}=
\begin{pmatrix}
	1 \\ 
	1 \\ 
	1\\ 
	\vdots\\ 
	1 \\ 
	1
\end{pmatrix}
,\;\;
{\rm and}\;\;
\mathbf{M}^{(k+1)\times (k+1)}=
\begin{pmatrix}
	1 & 1 & 0 &  0 &\cdots& 0 & 0 \\ 
	1 & 1 & 1 &  0 &\cdots&   0 & 0 \\ 
	1 & 1 & 1 &  1 &\cdots&  0 & 0\\ 
	1 & 1 & 1 &  1 &\cdots&  0& 0 \\ 
	\vdots&\vdots& \vdots&\vdots&\ddots& \vdots& \vdots\\ 
	1 & 1 & 1 & 1 &\cdots& 1& 1 \\ 
	1 & 1 & 1 & 1 &\cdots& 1 & 1
\end{pmatrix}.
$$

We know from \cite[Lemma 6]{NL_grow} that the characteristic polynomial of any recurrence sequence $(r)$ defined by the linear combination of the recurrence sequences $(a_{j})=(a_{i,j})_{i\geq0}$ of system \eqref{eqsystem}, moreover the characteristic polynomial of the coefficient matrix $\mathbf{M}$ of system \eqref{eqsystem} coincide. (The reader find a more precise theorem in \cite{NSz_Power}.) 
Consequently, we have to determine the characteristic polynomial $p_k(x)$ of $\mathbf{M}$, and then $p_k(x)$ yields the common recurrence relation of the sequences $\{a_{j}\}$ and their linear combinations.

Since $$p_k(x)=|x\mathbf{I}-\mathbf{M}|,$$
where $\mathbf{I}$ is the appropriate unit matrix,  
we obtain
\begin{eqnarray*}
	p_0(x)&=&x-1,\\
	p_1(x)&=&\begin{vmatrix} x-1 & -1 \\-1 & x-1\end{vmatrix}
	=x^2-2x,
\end{eqnarray*}	
and
\begin{eqnarray*}
	p_k(x)&=&
	\begin{vmatrix}
		x-1& -1 &  0 &  0 &\cdots& 0 & 0 \\ 
		-1 & x-1& -1 &  0 &\cdots&   0 & 0 \\ 
		-1 & -1 & x-1& -1 &\cdots&  0 & 0\\ 
		-1 & -1 & -1 & x-1&\cdots&  0& 0 \\ 
		\vdots&\vdots& \vdots&\vdots&\ddots& \vdots& \vdots\\ 
		-1 & -1 & -1 & -1 &\cdots& x-1& -1 \\ 
		-1 & -1 & -1 & -1 &\cdots& -1 & x-1
	\end{vmatrix}_{(k+1)\times(k+1)}\\
	&=&\allowdisplaybreaks
	\begin{vmatrix}
		x & -1 &  0 &  0 &\cdots& 0 & 0 \\ 
		-x &  x & -1 &  0 &\cdots&   0 & 0 \\ 
		0 & -x &  x & -1 &\cdots&  0 & 0\\ 
		0 &  0 & -x & x&\cdots&  0& 0 \\ 
		\vdots&\vdots& \vdots&\vdots&\ddots& \vdots& \vdots\\ 
		0 &  0 & 0 & 0 &\cdots& x& -1 \\ 
		0 &  0 & 0 & 0 &\cdots& -x & x-1
	\end{vmatrix}_{(k+1)\times(k+1)}\\
	&=&
	x\cdot
	\begin{vmatrix}
		x & -1 &  0 &  \cdots& 0 & 0 \\ 
		-x &  x & -1 &  \cdots&   0 & 0 \\ 
		0 & -x &  x & \cdots&  0 & 0\\ 
		\vdots&\vdots& \vdots&\ddots& \vdots& \vdots\\ 
		0 &  0 & 0 & \cdots& x& -1 \\ 
		0 &  0 & 0 & \cdots& -x & x-1
	\end{vmatrix}_{k\times k}
	+
	\begin{vmatrix}
		-x & -1 &  0 &  \cdots& 0 & 0 \\ 
		0 &  x & -1 &  \cdots&   0 & 0 \\ 
		0 & -x &  x & \cdots&  0 & 0\\ 
		\vdots&\vdots& \vdots&\ddots& \vdots& \vdots\\ 
		0 &  0 & 0 & \cdots& x& -1 \\ 
		0 &  0 & 0 & \cdots& -x & x-1
	\end{vmatrix}_{k \times k}\\
	&=&
	x\cdot
	p_{k-1}(x)
	-x\cdot
	\begin{vmatrix}
		x & -1 &   \cdots& 0 & 0 \\ 
		-x &  x &   \cdots&   0 & 0 \\ 
		\vdots& \vdots&\ddots& \vdots& \vdots\\ 
		0 &  0 &  \cdots& x& -1 \\ 
		0 &  0 &  \cdots& -x & x-1
	\end{vmatrix}_{(k-1)\times(k-1)}\\
	&=&x\cdot p_{k-1}(x) - x\cdot p_{k-2}(x).
\end{eqnarray*}
In the calculation above, first we subtracted the second column of the determinant from the first column, then we subtracted the third column from the second one, and so on. Secondly, we expanded the determinant by its first row. Thirdly, we expanded the second determinant by its first column. Finally, for $p_k(x)$  we have the binary recurrence relation 
\begin{equation}\label{eq:rec_kx}
	p_k(x)= x\cdot p_{k-1}(x) - x\cdot p_{k-2}(x).
\end{equation}

To derive the explicit form of $p_k(x)$ we solve the characteristic equation $z^2-xz+x=0$ of \eqref{eq:rec_kx}. Note that the left-hand side is a quadratic polynomial in $z$.  The solutions are 
$z_{1}=\frac{1}{2}(x+\sqrt{x(x-4)})$ and $z_{2}=\frac{1}{2}(x-\sqrt{x(x-4)})$.

By the fundamental theorem of homogenous linear recurrences $p_k(x)=\alpha z_1^k + \beta z_2^k$, where $\alpha $ and $\beta$ are determined by the linear equation system  
\begin{eqnarray*}
	p_0(x) &=& \alpha  + \beta, \\ 
	p_1(x)  &=& \alpha z_1 + \beta z_2. 
\end{eqnarray*}
Obviously,
\begin{equation}\label{alfabeta}
	\begin{aligned}
		\alpha&= \frac12 \,{\frac {{x}^{2}-5x+4+ (x-3) \sqrt{x(x-4)}}{x-4}},\\
		\beta&= \frac12 \,{\frac {{x}^{2}-5x+4 -(x-3) \sqrt{x(x-4)}}{x-4}}.
	\end{aligned}
\end{equation}

Hence we have proved the following theorem.
\begin{theorem}
	The characteristic polynomial $p_k(x)$ 
	has degree $k+1$, it satisfies \eqref{eq:rec_kx}, its explicit formula with \eqref{alfabeta} is
	\begin{equation*}
		p_k(x)=\alpha \left(\frac{1}{2}(x+\sqrt{x(x-4)})\right)^k + \beta  \left(\frac{1}{2}(x-\sqrt{x(x-4)})\right)^k.
	\end{equation*} 
\end{theorem}

Because each recurrence coefficient in \eqref{eq:rec_kx} is one of  $\pm x$, the factorization of $p_k(x)$ contains a factor $x^m$ for some positive integer $m$. The next theorem provides, among others, the precise exponent $m$ in the factorization of $p_k(x)$. 

\begin{theorem} \label{th:charpoly}
	The characteristic polynomials $p_k(x)$ can be given by 
	\begin{equation*}  
		p_k(x)=x^{\Cl{k}} \sum_{i=0}^{\Fl{k}}  (-1)^i \binom{k+2-i}{i}x^{\Fl{k}-i}, \qquad k\geq 0.
	\end{equation*}
\end{theorem}

\begin{proof}
	Observe that $\left\lceil \frac{k}{2}\right\rceil+\left\lfloor \frac{k}{2}\right\rfloor=k$, and $\left\lfloor \frac{k}{2}\right\rfloor+1=\left\lfloor \frac{k+2}{2}\right\rfloor$. Put $n:=k+2$. Thus,
	\begin{equation*}
		p_{n-2}(x)=\sum_{i=0}^{\left\lfloor n/2\right\rfloor}  \binom{n-i}{i}(-1)^ix^{n-1-i}= \frac{1}{x} \sum_{i=0}^{\left\lfloor n/2\right\rfloor}  \binom{n-i}{i}x^{n-2i}(-x)^i.
	\end{equation*}

	Now we shall apply Theorem 1 of \cite{BKSz}, but first, for the convenience of the readers we present it in
	\begin{lemma}[{\cite[Theorem 1]{BKSz}}]\label{lem:BKSz}
		The terms of the sequence $(T_n)_n$ given by 
		$$ T_{n+1}=  \sum_{i=0}^{\left\lfloor (n-p)/(q+r)\right\rfloor}  \binom{n-qi}{p+ri}x^{n-p-(q+r)i}y^{p+ri}$$
		satisfy the linear recurrence relation
		$$T_n-x\binom{r}{1}T_{n-1}+x^2\binom{r}{2}T_{n-2}+\cdots+ (-1)^rx^r\binom{r}{r}T_{n-r}=y^rT_{n-r-q}.$$
	\end{lemma}
	
	We choose the parameters of Lemma \ref{lem:BKSz} as $r=q=1$, $p=0$, and $y=-x$. Hence the terms $T_{n+1}=\sum_{i=0}^{\left\lfloor n/2\right\rfloor}  \binom{n-i}{i}(-x)^ix^{n-2i}$ satisfy $T_{n}=xT_{n-1}+(-x)T_{n-2}$ for $n\ge3$. 
	Consequently, $p_k(x)=T_{n-2}(x)/x$ implies
	\begin{equation} \label{beszurt}
		p_{k+2}(x)=xp_{k+1}-xp_{k},\qquad k\ge 1.
	\end{equation}

	Since the initial polynomials are 
	$$
	p_0(x)=\frac{1}{x}T_3(x)=\frac{1}{x}(x^2-x)=x-1,\quad{\rm and}\quad p_1(x)=\frac{1}{x}T_4(x)=\frac{1}{x}(x^3-2x^2)=x^2-2x,
	$$
	then (\ref{beszurt}) is true for $k=0$, too.
\end{proof}

For example, the first few characteristic polynomials are
\begin{equation*}
	\begin{aligned}
		p_0(x)&=x-1,\\
		p_1(x)&=x(x-2),\\
		p_2(x)&=x (x^2-3 x+1),\\
		p_3(x)&=x^2(x^2-4x+3),\\
		p_4(x)&=x^2(x^3-5x^2+6x-1),\\
		p_5(x)&=x^3(x^3-6x^2+10x-4),\\
		p_6(x)&=x^3(x^4-7x^3+15x^2-10x+1),\\
		p_7(x)&=x^4(x^4-8x^3+21x^2-20x+5),\\
		p_8(x)&=x^4(x^5-9x^4+28x^5-35x^2+15x-1), \\
		p_9(x)&=x^5(x^{5}-10x^4+36x^3-56x^2+35x-6),\\
		p_{10}(x)&=x^5(x^{6}-11x^5+45x^4-84x^3+70x^2-21x+1).
	\end{aligned}
\end{equation*}

\section{Proofs of the main theorems}

\subsection{Right-down diagonal sequences}

Now we are ready to give the recurrences of the right-down diagonal sequences $A_j^{(k)}$ by the help of Theorem~\ref{th:charpoly}.
\begin{theorem} 
	In case of a fixed $j$ we have the recurrence relation
	\begin{equation}\label{eq:rec_a_ij}
		0= \sum_{i=0}^{\Fl{k}}  (-1)^i \binom{k+2-i}{i} a_{n-i,j}, \qquad n\geq \Fl{k}, 0\leq j\leq k+1.
	\end{equation}
\end{theorem}
\begin{proof} 
	Recall that $p_k(x)$ is not only the characteristic polynomial of $\mathbf{M}$  but also of any recurrence sequence of system \eqref{eqsystem}, moreover,  of any  linear combination sequence of them. 
	This is why, if we substitute the power $x^i$ by $a_{i,j}$ in $p_k(x)=0$ in Theorem~\ref{th:charpoly},  then since $x\ne 0$, we get
	\begin{equation*}
		\begin{aligned}
			0&= \sum_{i=0}^{\Fl{k}}  (-1)^i \binom{k+2-i}{i} a_{{\Fl{k}-i},j}
			= \sum_{i=0}^{\Fl{k}}  (-1)^i \binom{k+2-i}{i} a_{{\Fl{k}-i}+n-\left({\Fl{k}}\right),j}\\
			&= \sum_{i=0}^{\Fl{k}}  (-1)^i \binom{k+2-i}{i} a_{n-i,j}.
		\end{aligned}
	\end{equation*}
\end{proof}

Expressing the item $a_{n,j}$ from recurrence equation \eqref{eq:rec_a_ij} we obtain the result of Theorem~\ref{th:main1}.
For example, the first few recurrence relations are
\begin{equation*}
	\begin{aligned}
		k=0:&\quad&a_{n,j}&=a_{n-1,j},\\
		k=1:& &a_{n,j}&=2a_{n-1,j},\\
		k=2:& &a_{n,j}&=3a_{n-1,j}-a_{n-2,j},\\
		k=3:& &a_{n,j}&=4a_{n-1,j}-3a_{n-2,j},\\
		k=4:& &a_{n,j}&=5a_{n-1,j}-6a_{n-2,j}+a_{n-2,j},\\
		k=5:& &a_{n,j}&=6a_{n-1,j}-10a_{n-2,j}+4a_{n-2,j},\\
		k=6:& &a_{n,j}&=7a_{n-1,j}-15a_{n-2,j}+10a_{n-2,j}-a_{n-3,j},\\
		k=7:& &a_{n,j}&=8a_{n-1,j}-21a_{n-2,j}+20a_{n-2,j}-5a_{n-3,j},\\
		k=8:& & a_{n,j}&=9a_{n-1,j}-28a_{n-2,j}+35a_{n-2,j}-15a_{n-3,j}+a_{n-3,j},\\ 
		k=9:& & a_{n,j}&=10a_{n-1,j}-36a_{n-2,j}+56a_{n-2,j}-35a_{n-3,j}+6a_{n-4,j},\\
		k=10:& & a_{n,j}&=11a_{n-1,j}-45a_{n-2,j}+84a_{n-2,j}-70a_{n-3,j}+21a_{n-4,j}-a_{n-5,j}.
	\end{aligned}
\end{equation*}

\subsection{Zig-zag sequences}

Theorem~\ref{th:main2} is the simple corollary of Theorem~\ref{th:main1}. Indeed, we only merge two sequences satisfying the same recurrence relation. Hence the statement is obvious.

\section{Sum of rows, columns, and left-down diagonal sequences}

Let $R^{(k)}=(r^{(k)}_n)$ be the sum sequence of the values of the $n$th row of square $k$--zig-zag shape. Considering the partial sum relation~\eqref{eq:partsum} we obtain
$$r^{(k)}_n= \sum_{j=0}^{k+1} a_{n,j}=a_{n,0}+a_{n+1,k}.$$ 
So, the recurrence sequence $r^{(k)}_n$ is the linear combination of sequences $A^{(k)}_n$, therefore they have the same characteristic polynomial and the same recurrence relation.  

Let $C^{(k)}=(c^{(k)}_n)$ be the sum sequence of columns.
As $a_{n+1,0}=a_{n,1}=a_{n,0}+a_{n-1,2}=a_{n,0}+a_{n-1,1}+a_{n-2,3}=\cdots=\sum_{j=0}^{\min\{n,k+1\}-1}a_{n-j,j}$, then 
\begin{equation*}
	c^{(k)}_n= \sum_{j=0}^{\min\{n,k+1\}} a_{n-j,j}=\begin{cases}
		a_{n+1,0}, & \text{ if } n\leq k; \\ 
		a_{n+1,0}+a_{n-k-1,k+1}, & \text{ if } n>k. 
	\end{cases}
\end{equation*}

Let $D^{(k)}=(d^{(k)}_n)$  be the left-down diagonal sequence,  where
\begin{equation*}\label{eq:w}
	d_{n}=\begin{cases}
		\ \displaystyle\sum_{\ell=0}^{\ell\leq \frac{n}{2},\,  2\ell\leq k}a_{\frac{n}{2}-\ell,2\ell}, & \text{ if $n$ is even}; \\[1.5em]
		\displaystyle\sum_{\ell=0}^{\ell\leq \frac{n}{2},\,  2\ell\leq k+1}a_{\frac{n}{2}-\ell,2\ell+1}, & \text{ if $n$ is odd}. 
	\end{cases}
\end{equation*}

Since all the $A_j^{(k)}$ sequence satisfy the same recurrence relation, then $C^{(k)}$ and $D^{(k)}$  are so.

\begin{figure}[H]
	\centering
	\includegraphics[scale=0.5]{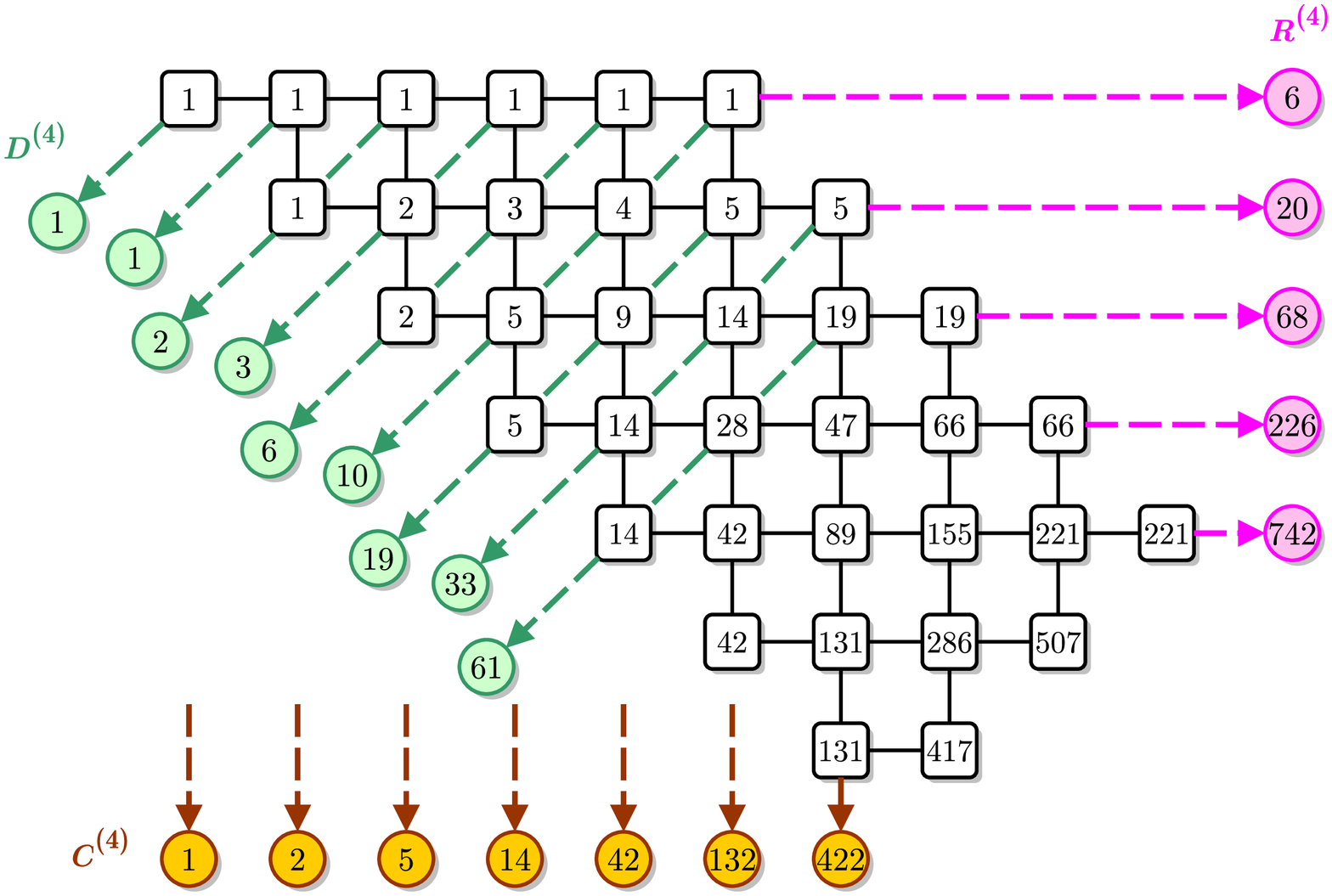} 
	\caption{Sum of rows, columns, and left-down diagonal sequences}
	\label{fig:zigzag_k4_diag}
\end{figure}

For example, see Figure~\ref{fig:zigzag_k4_diag} when $k=4$. 
\begin{equation*}
	\begin{aligned}
		R^{(4)} &= (6, 20, 68, 226, 742, 2422, 7884, 25630, 83268, 270444,  \ldots) &=& \text{ not in the OEIS}, \\
		C^{(4)} &= (1,2,5,14,42,132,422,1360,4400,14262,46270,  \ldots) &=& \text{ not in the OEIS}, \\
		(d^{(k)}_{2n}) &= (1, 2, 6, 19, 61, 197, 638, 2069, 6714, 21794, 70755,  \ldots) &=&\ \text{\seqnum{A052975}},\\
		(d^{(k)}_{2n+1}) &= (1, 3, 10, 33, 108, 352, 1145, 3721, 12087, 39254, 127469, \ldots) &=&\  \text{\seqnum{A060557}},\\
		D^{(4)} &= (1, 1, 2, 3, 6, 10, 19, 33, 61, 108, 197, 352, 638, 1145, 2069,  \ldots) &=&\ \text{\seqnum{A028495}}.
	\end{aligned}
\end{equation*}

\section{Examples}\label{sec:example}

In this section, we give some example sequences with different $k$ values appearing in the {\it On-Line Encyclopedia of Integer Sequences} \cite{OEIS} (see Table~\ref{tab:ex1} and Table~\ref{tab:ex2}).
\begin{table}[H]\centering
	\ra{1.3}
	\begin{tabular}{@{}rlclclcl@{}}\toprule
		$k=0$ 
		& \multicolumn{6}{l}{$A^{(0)}_0=A^{(0)}_1$  =\seqnum{A000012}.}  \\
		$k=1$ 
		& $A^{(1)}_1$=\seqnum{A000079}\\
		$k=2$ 
		& $A^{(2)}_0$=\seqnum{A001519} & & $A^{(2)}_1 $=\seqnum{A001519}$^*$&& $A^{(2)}_2$=\seqnum{A001906}$^*$ \\
		$k=3$ 
		& $A^{(3)}_0$=\seqnum{A124302} & & $A^{(3)}_1 $=\seqnum{A000244} && $A^{(3)}_2$=\seqnum{A000244}&& $A^{(3)}_3$=\seqnum{A003462}\\
		$k=4$ 
		& \multicolumn{6}{l}{These sequences are presented in Section~\ref{sec:2}.}  \\
		$k=5$
		& $A^{(5)}_0$=\seqnum{A024175} & & $A^{(5)}_1$=\seqnum{A024175}$^*$ && $A^{(5)}_2$=\seqnum{A094803}&& $A^{(5)}_3=$\seqnum{A007070}\\
		& $A^{(5)}_4$=\seqnum{A094806} & & $A^{(5)}_5$=\seqnum{A094811} \\
		$k=6$
		& $A^{(6)}_0$=\seqnum{A080938} & & $A^{(6)}_1$=\seqnum{A080938}$^*$ && $A^{(6)}_2$=\seqnum{A094826}   & & $A^{(6)}_3=$\seqnum{A094827}\\
		& $A^{(6)}_4$=\seqnum{A094828} & & $A^{(6)}_5$=\seqnum{A094829} & & $A^{(6)}_6$=\seqnum{A094256}\\
		$k=7$
		& $A^{(7)}_0$=\seqnum{A033191} & & $A^{(7)}_1$=\seqnum{A033191}$^*$ & &   	$A^{(7)}_2$=\seqnum{A033190}   & & $A^{(7)}_3=$\seqnum{A094667}\\
		& \begin{tabular}{@{}c@{}}$A^{(7)}_4$=\seqnum{A030191}; \\ \hphantom{$A^{(7)}_4$=}\seqnum{A093131}$^*$ \end{tabular}   & & $A^{(7)}_5$=\seqnum{A094788} & & $A^{(7)}_6$=\seqnum{A094825} & & $A^{(7)}_7$=\seqnum{A094865}$^{**}$ \\
		$k=8$
		& $A^{(8)}_0$=\seqnum{A211216}  & & $A^{(8)}_1$=\seqnum{A211216}$^*$ & &   	$A^{(8)}_4$=\seqnum{A224422}   & & $A^{(8)}_5=$\seqnum{A221863}\\
		$k=9$ & $A^{(9)}_5$=\seqnum{A216263}  \\ \bottomrule
		\multicolumn{6}{l}{$^* n\geq1$, $^{**} n\geq3$.}
	\end{tabular}
	\caption{Sequences $A^{(k)}_j$ appearing in the OEIS} 
	\label{tab:ex1}
\end{table}

\begin{table}[H]\centering
	\ra{1.3}
	\begin{tabular}{@{}rlclclcl@{}}\toprule
		$k=0$ 
		& \multicolumn{6}{l}{$Z^{(0)}_0$=\seqnum{A000001}}  \\
		$k=1$ 
		& \multicolumn{4}{l}{$Z^{(2)}_1$=\seqnum{A000045}, Fibonacci sequence}
		&$Z^{(2)}_2$=\seqnum{A001906}$^*$  \\
		$k=3$ 
		& $Z^{(3)}_0$=\seqnum{A124302} & & $Z^{(3)}_1$=\seqnum{A000244} && $Z^{(3)}_2$=\seqnum{A232801}$^*$  & & $Z^{(3)}_3$=\seqnum{A052993}\\
		$k=4$ 
		& \multicolumn{6}{l}{These sequences are presented in Section~\ref{sec:2}.}  \\
		\bottomrule
		\multicolumn{6}{l}{$^* n\geq1$.}
	\end{tabular}
	\caption{Sequences $Z^{(k)}_j$ appearing in the OEIS} 
	\label{tab:ex2}
\end{table}

\section{Acknowledgments}

For L. Sz., {this work was supported by Hungarian National Foundation for Scientific Research Grant Nos.~128088 and 130909.}
The authors would like to thank the anonymous referee for carefully reading the manuscript, and for the useful suggestions and improvements.

\end{document}